\documentclass[11pt]{article}
\usepackage[ansinew]{inputenc}
\usepackage{latexsym}
\usepackage{textcomp}
\usepackage{color}
\usepackage{amscd}

\input{amssym.def}
\input{amssym}

     \newcommand{\fg}{\goth{g}}
     
     \newcommand{\fr}{\goth{r}}

     \newcommand{\E}{\Bbb{E}}
     \newcommand{\F}{\Bbb{F}}
     
     \newcommand{\Ll}{\Bbb{L}}  

     \newcommand{\Z}{\Bbb{Z}}

    \newcommand{\ol}[1]{\overline{#1}}

    \newcommand{\ti}[1]{\tilde{#1}}
    \newcommand{\group}[1]{\langle{#1}\rangle}
    
    \newcommand{\op}{\oplus}
    \newcommand{\ot}{\otimes}
    \newcommand{\me}{^{-1}}
    \newcommand{\mal}{^{\times}}
    
    \newcommand{\df}{\stackrel{\mathrm{def}}{=}}
    
    \newcommand{\mr}{\mathrm}
    
    \newcommand{\ilim}[1]{\raisebox{-1mm}{$\lim\atop{\leftarrow\atop\raisebox{0.5mm}{{\tiny $#1$}}}$}}

    \newcommand{\zl}{{\Bbb{Z}_l}}
    
    \newcommand{\ql}{{\Bbb{Q}_l}}

    \newcommand{\into}{\rightarrowtail}
    \newcommand{\onto}{\twoheadrightarrow}
    \newcommand{\lto}{\longrightarrow}
    \newcommand{\da}{\downarrow}

    \newcommand{\pht}{\phantom}

    \def\daz#1{#1\da\pht{#1}}
    
    \newcommand{\Ga}{\Gamma}

    \newcommand{\De}{\Delta}
    \newcommand{\ze}{\zeta}
    \newcommand{\la}{\lambda}
    \newcommand{\La}{\Lambda}

    \newcommand{\sn}{\par\smallskip\noindent}
    \newcommand{\mn}{\par\medskip\noindent}
    \newcommand{\bn}{\par\bigskip\noindent}
    \newcommand{\bbn}{\par\bigskip\bigskip\noindent}
    
    \newcommand{\Section}[2]{\bbn {\large #1\,. \ {\sc #2}}
                             \nopagebreak
                             \nz}

    \newcommand{\nf}[2]{\\[1.5ex]
                        \bmp{1cm}
                         (#1)
                        \emp 
                        \bmp{13.5cm}
                         \bct
                          $#2$
                         \ect
                        \emp\\[1.5ex]
         }

    \newcommand{\sss}{\scriptstyle}
    \newcommand{\nz}{\\[1ex]}

    \newcommand{\hsp}[1]{\hspace*{#1mm}}

    \newcommand{\mmargin}{
     \textheight 230truemm
     \textwidth 155truemm
     \topmargin -10truemm
     \oddsidemargin 5truemm
     \evensidemargin 5truemm
     }

    \newcommand{\bmp}{\begin{minipage}}
    \newcommand{\emp}{\end{minipage}}
    \newcommand{\btb}{\begin{tabular}}
    \newcommand{\etb}{\end{tabular}}
    \newcommand{\barr}{\begin{array}}
    \newcommand{\earr}{\end{array}}
    \newcommand{\bit}{\begin{itemize}}
    \newcommand{\eit}{\end{itemize}}
    \newcommand{\ben}{\begin{enumerate}}
    \newcommand{\een}{\end{enumerate}}
    \newcommand{\bct}{\begin{center}}
    \newcommand{\ect}{\end{center}}
    \newcommand{\bfr}{\begin{flushright}}
    \newcommand{\efr}{\end{flushright}}
    \newcommand{\bea}{\begin{eqnarray*}}
    \newcommand{\eea}{\end{eqnarray*}}
    \newcommand{\bqo}{\begin{quote}}
    \newcommand{\eqo}{\end{quote}}
    \newcommand{\bdc}{\begin{description}}
    \newcommand{\edc}{\end{description}}
    \newcommand{\bdia}{\begin{CD}}
    \newcommand{\edia}{\end{CD}}

    \definecolor{light}{gray}{.3}

    \newcommand{\coker}{\mathrm{coker\,}}

    \newcommand{\im}{\mathrm{im\,}}
    \newcommand{\mod}{\mathrm{\ mod\ }}

    \newcommand{\sr}[2]{{\,\stackrel{#1}{#2}\,}}
    \newcommand{\fra}[2]{{\,\frac{#1}{#2}\,}}

    \newcommand{\Stop}{\edc \sn\rm} 
    \newcommand{\Stopp}{\Stop\vspace*{-10mm}}
    
    \newcommand{\Lemma}[1]{\sn
           \bdc
           \item[{\sc Lemma {#1}.}] \em }
    
    \newcommand{\Proposition}[1]{\sn
           \bdc
           \item[{\sc Proposition {#1}.}] \em }
    \newcommand{\corollary}{\sn
           \bdc
           \item[{\sc Corollary.}] \em }
    \newcommand{\Corollary}[1]{\sn
           \bdc
           \item[{\sc Corollary {#1}.}] \em }

    \newcommand{\remark}{\sn{\sc Remark.} \ }
    
    \newcommand{\definition}{\sn
           \bdc
           \item[{\sc Definition.}] \em }

\mmargin

\def\ab{{^\mr{ab}}}
\def\law{{\La_\wedge}}
\def\lawm{{\La_\wedge\mal}}
\def\lab{{\La_\bullet}}
\def\zlt{{\zl[[T]]}}
\def\rad{{\mr{rad}}}
\def\lawh{{\La_\wedge[H]}}
\def\up{{\ol 1-\ol\psi}}

\begin{document}

\title{The integral logarithm in Iwasawa theory: an exercise}
 \author{Jürgen Ritter \ $\cdot$ \ Alfred Weiss \
 \thanks{We acknowledge financial support provided by NSERC and the University of Augsburg.}
 }
\date{\pht\today}

\maketitle




 \bbn Let $\La=\zl[[T]]$ denote the ring of power series in one
 variable over the $l$-adic integers $\zl$, where $l$ is an odd
 prime number. We localize $\La$ at the prime ideal $l\cdot\La$ to
 arrive at $\La_\bullet$ and then form the completion
 $$\La_\wedge=\ilim{n}\La_\bullet/l^n\La_\bullet\ .$$ The integral
 logarithm \ $\Ll:\La_\wedge\mal\to\La_\wedge$ \ is defined by
 $$\Ll(e)=\fra1l\log\fra{e^l}{\psi(e)}\,,$$ where
 $\psi:\La_\wedge\to\La_\wedge$ is the $\zl$-algebra homomorphism induced by
 $\psi(T)=(1+T)^l-1$ and with `log' defined by the usual power
 series.
 \sn In this paper, the unit group $\La_\wedge\mal$ as well as
 $\ker(\Ll)$ and $\coker(\Ll)$ are studied -- more precisely, we
 study the analogous objects when $\La_\wedge$ is replaced by the
 group ring $\La_\wedge[H]$ of a finite abelian $l$-group $H$.
 \sn The interest in doing so comes from recent work in Iwasawa theory in
 which refined `main conjectures' are formulated in terms of
 the $K$-theory of completed group algebras $\zl[[G]]$ with $G$ an
 $l$-adic Lie group (see [FK], [RW2]). For $l$-adic Lie groups of
 dimension 1, use of the integral logarithm $\Ll$ has reduced the
 `main conjecture' to questions of the existence of special
 elements (``pseudomeasures'') in $K_1(\zl[[G]]_\bullet)$\,, by
 Theorem A of [RW3], and, more recently ([RW4,5]), to still open
 {\em logarithmic congruences} between Iwasawa $L$-functions.
 Moreover, $\Ll$ has been indispensable for the proof of the `main
 conjecture' in the few special cases ([K]\,,\,[RW6]) which have
 been settled so far.
 \sn Apparently the integral
 logarithm $\Ll$ when applied to $K_1(\zl[[T]]_\bullet)$ takes its
 values only
 in $\zl[[T]]_\wedge$, and here it helps, as it did for
 finite $G$ (see [O] and [F]), to obtain structural information
 about $K_1$, which, in particular, implies that $\coker(\Ll)$ can
 be detected on the abelianization $G\ab$ of $G$ (by [RW3, Theorem
 8]). Then \bit\item[]$G\ab=H\times\Ga$, with $\Ga\simeq\zl$ and $H$ as
 before, \item[] $K_1(\zl[[G\ab]]_\wedge)=\zl[[G\ab]]_\wedge\mal$ \quad(see [CR, 40.31 and 40.32 (ii)])\,, \item[]
 $\zl[[G\ab]],\,\zl[[G\ab]]_\wedge$ are $\La[H]$ and $\La_\wedge[H]$,
 respectively, \item[] and $\psi$ is induced by the map $g\mapsto g^l$ on
 $G\ab$.\eit
 For these reasons it seems worthwhile to present a rather complete
 understanding of $\Ll$ in the abelian situation, which is the
 purpose of our exercise.
 \bn The content of the paper is as follows. In a first section we
 consider $\La$ and define an integral exponential
 $\E$ on $T^2\La$ which is inverse to $\Ll$ (on $1+T^2\La$). As a
 consequence, we obtain the decomposition
 $$\La\mal=\zl\mal\times(1+T)^\zl\times\E(T^2\La)$$ for the unit
 group $\La\mal$ of $\La$ (which reminds us of [C, Theorem 1]).
 Applying  $\Ll$ to the decomposition yields a generalization of
 the Oliver congruences [O, Theorem 6.6].
 \sn The second section centers around $\La_\wedge$ and two
 important subgroups $$\Xi=\{\sum_{k=-\infty}^\infty
 x_kT^k\in\law:x_k=0\ \mr{for}\ l|k\}\quad \mr{and}\quad
 \Xi_2=\{\sum_{k\ge 2}x_kT^k\in\Xi\}\ .$$ In terms of these we
 exhibit natural decompositions of $\law$ and $\lawm$, which
 leads immediately to $\ker(\Ll)$ and $\im(\Ll)\,,\,\coker(\Ll)$.
 \sn Section 3 is still concerned with $\law$\,: we determine the
 kernel and cokernel of its endomorphism $1-\psi$.
 \sn This will be used in the last section, \S4, where we extend most of the results to
 the group ring $\law[H]$ of a finite abelian $l$-group $H$ over
 $\law$ and determine $\ker(\Ll)$ and $\coker(\Ll)$ here.
 \Section{1}{The integral exponential $\E$ and $\La\mal$}
 Recall that $\La$ is the ring $\zlt$ of formal power series
 $\sum_{k\ge0}y_kT^k$ with coefficients $y_k\in\zl$\,, and that the
 integral logarithm is defined on the units $e\in\La\mal$ of $\La$ by
 $$\Ll(e)=\fra1l\log\fra{e^l}{\psi(e)}\quad\mr{where}\quad\psi(T)=(1+T)^l-1\
 .$$ Moreover, note that $1+T^2\La$ is a subgroup of $\La\mal$ since
 $T\in\rad(\La)=\group{l,T}$\,.
 \sn We now turn to the integral exponential $\E$ on $T^2\La$\,: This is the formal
 power series, with coefficients in $\ql$, defined by
 $$\E(y)=\exp\Big(\sum_{i\ge0}\fra{\psi^i(y)}{l^i}\Big)\in\ql[[T]]\quad\mr{for\ each}\quad y\in T^2\La\ .$$
 Observe that $\E$ and $\psi$ commute.
 \Lemma{1} $\E(y)\in 1+T^2\La$\,, and $\E$ and $\Ll$ are inverse to
 each other $$T^2\La \ {{\E\atop\rightleftarrows}\atop\sss{\Ll}} \ 1+T^2\La\ .$$\Stop
 The proof is an adaptation of that of the
 Dwork-Dieudonn\'e lemma (see [L, 14\,\S2])\,: if
 $f(T)\in1+T^2\ql[[T]]$ satisfies
 $\fra{f(T)^l}{\psi(f(T))}\in1+lT^2\zlt$\,, then $f(T)\in
 1+T^2\zlt$\,.
 \sn First, $\psi^i(T)\equiv l^iT\mod T^2\zlt$ implies
 $\psi^i(T^k)\equiv l^{ik}T^k\mod T^{k+i}\zlt$\,, and thus, if
 $y\in y_kT^k+T^{k+1}\zlt$ with $y_k\in\zl$, then $$\barr{l}\psi^i(y)/l^i\in
 y_kl^{(k-1)i}T^k+T^{k+1}\ql[[T]]\,;\quad\mr{whence,\ for}\ k\ge2\,,\\[1.5mm]
 \sum_{i\ge0}\fra{\psi^i(y)}{l^i}\in\Big(\sum_{i\ge0}y_kl^{(k-1)i}T^k\Big)+T^{k+1}\ql[[T]]=
 y_k(1-l^{k-1})\me T^k+T^{k+1}\ql[[T]]\,,\\[1.5mm]
 \mr{so}\quad \E(y)\in 1+y_k(1-l^{k-1})\me T^k+T^{k+1}\ql[[T]]\ .\earr$$
 Second,
 \nf{$\ast$}{\E(y)^l/\E(\psi(y))=\exp\Big(l\sum_{i\ge0}\fra{\psi^i(y)}{l^i}-\sum_{i\ge0}\fra{\psi^{i+1}(y)}{l^i}\Big)=
 \exp(ly)\in 1+lT^2\zlt\,,} which brings us in a position to
 employ the Dwork-Dieudonn\' e argument to obtain
 $\E(y)\in1+T^2\zlt=1+T^2\La$ and, in particular, $\E(y)\in 1+y_k(1-l^{k-1})\me T^k+T^{k+1}\zlt$\,.
 \sn Moreover, given $b_k\in\zl$ and setting
 $a_k=(1-l^{k-1})b_k$, then
 $\E(a_kT^k)\in1+b_kT^k+T^{k+1}\zlt$ for all $k$, which implies that
 $\E(T^2\La)=1+T^2\La$\,.
 \sn We finish the proof of the lemma by showing $\Ll\E(y)=y\,,$
 and $\E\Ll(1+y)=1+y$ whenever $y\in T^2\La$\,, so $1+y=\E(\ti y)$\,:
 $$\barr{l}
 \Ll\E(y)=\fra1l\log\fra{\E(y)^l}{\E(\psi(y))}\sr{(\ast)}{=}\fra1l\log\exp(ly)=y\,,\\[1.2mm]
 \E\Ll(1+y)=\E\Ll(\E(\ti y))=\E(\ti y)=y\,.\earr$$
 \Corollary{1} \
 $\barr{l} \Ll(1+l\La)=(l-\psi)\La\\
 \exp(ly)=(1+T)^{ly_1}\E((l-\psi)y) \  {i\!f} \ y\equiv y_1T\mod
 T^2\La \ (and\ y_1\in\zl)\earr$\Stop
 Since $\exp(l\La)=1+l\La$, for the first assertion it suffices to
 compute $$\Ll(\exp(ly))=
 \fra1l\log(\exp(ly)^{l-\psi})=\fra1l\log\exp(l(l-\psi)y)=(l-\psi)y\,.$$
 The second assertion holds for $y=T$ (so $y_1=1)$\,:
 $$(1+T)^{-l}\exp(lT)\in1+T^2\zlt\,,\ \mr{hence}\ (1+T)^{-l}\exp(lT)=\E(z)\
 \mr{for\ some}\ z\in T^2\zlt\,.$$
 Apply $\Ll$ and get $(l-\psi)(T)=z$ from the last but one displayed formula and
 as
 $\Ll(1+T)=\fra1l\log\fra{(1+T)^l}{1+(1+T)^l-1}=0$\,.
 \sn Next, take $y\in T^2\zlt$, so $\exp(ly)\in 1+T^2\zlt$ and
 again $\exp(ly)=\E((l-\psi)y)$.
 \sn The two special cases can be combined on writing
 $y=y_1T+(y-y_1T)$. This finishes the proof of the corollary.
 \mn Denote by $\mu_{l-1}$ the group of roots of unity in $\zl\mal$.
 \Corollary{2}  \
 $\barr{l}\La\mal=\zl\mal\times(1+T)^\zl\times\E(T^2\La)\\
 \ker(\Ll)=\mu_{l-1}\times (1+T)^\zl\,, \  \im(\Ll)=\zl\op T^2\La\earr$\Stop
 The first coefficient $e_0$ of
 $e=\sum_{k\ge0}e_kT^k\in\La\mal$ is a unit in $\zl$. Replacing
 $e$ by $e_0\me\cdot e=1+e_1'T+\cdots$ and then multiplying by
 $(1+T)^{-e_1'}$ gives the new unit $1+\ti
 e_2T^2+\cdots\in1+T^2\La$\,. Thus, by Lemma 1,
 $\La\mal=\zl\mal\cdot(1+T)^\zl\cdot\E(T^2\La)$\,, and the product
 is obviously direct. Since, on $\zl\mal$, $\Ll(\ze)=0$ precisely for
 $\ze\in\mu_{l-1}$, and since
 $\Ll(1+T)=0$\,, we also
 get the claimed description of the kernel and image of $\Ll$.
 \Section{2}{The integral logarithm $\Ll$ on $\law$}
 We recall that $\lab$ denotes the localization of $\La$ at the
 prime ideal $l\La$ and that $\law=\ilim{n}\lab/l^n\lab$\,. In
 particular, $\lab$ and $\law$ have the same residue field
 $\F_l((\ol T))$ (which carries the natural $\ol T$-valuation
 $v_{\ol T}$). It follows that
 $$\law=\{x=\sum_{k\in\Z}x_kT^k:x_k\in\zl\,,\
 \lim_{k\to-\infty}x_k=0\}\,.$$
 Such large things are basic objects in the theory of higher
 dimensional local fields [FeKu]; the map $\psi$ on $\law$ is extra
 structure which remembers the group $\Ga$.
 \sn In what follows we frequently use the decomposition
 $$\law=\law^-\op\zl\op\law^+\,,\quad\mr{where}\
 \law^\pm=\left\{\barr{l}\{x\in\law:x_k=0\ \mr{for}\ k\le0\}\\
 \{x\in\law:x_k=0\ \mr{for}\ k\ge0\}\,.\earr\right.$$
 Note that the three summands are subrings which are preserved by
 $\psi$. As a consequence, we see that
  \ $\La\cap(l-\psi)\law=(l-\psi)\La$\,.
 \definition $\barr{l}\Xi=\{x=\sum_{k\in\Z}x_kT^k\in\law:x_k=0\ \mr{when}\
 l\ \mr{divides}\ k\}\,,\\ \Xi_s=\{x=\sum_{k\ge
 s}x_kT^k\in\Xi\}\,,\quad\mr{where}\
 s\in\Z\,.\earr$\Stop
 \Lemma{2} \ben\item $l-\psi$ is injective on $\law$ and has image
 $\Ll(1+l\law)$\,,
 \item $\law=\Xi\op(l-\psi)\law$\een\Stop
 For the first assertion we make use of the commuting diagram
 \nf{1}{\barr{ccccc}\law&\sr{l}{\into}&\law&\onto&\F_l((\ol T))\\
 \daz{l-\psi}&&\daz{l-\psi}&&\daz{-\ol\psi}\\
 \law&\sr{l}{\into}&\law&\onto&\F_l((\ol T))\earr}
 with exact rows and with $\ol\psi(\ol T)=\ol T^l$\,, so
 $\ol\psi(\ol x)=\ol x^l$ for $\ol x\in\F_l((\ol T))\,. $ In
 particular, $-\ol\psi$ is injective and hence the snake lemma
 implies $\ker(l-\psi)=l\cdot\ker(l-\psi)$ from which
 $\ker(l-\psi)=0$ follows by $\bigcap_{n\ge 0}l^n\law=0$\,.
 \sn Regarding the image, we observe that $\exp(l\law)=1+l\law$ and
 recall $\Ll(\exp(ly))=(l-\psi)y$ from the proof of Corollary 1 of Lemma 1 (but now with $y\in\law$).
 \sn For the second assertion we make use of the commuting diagram
 \nf{2}{\barr{ccccc}\Xi&\sr{l}{\into}&\Xi&\onto&\Xi/l\Xi\\
 \da&&\da&&\da\\
 \law/(l-\psi)\law&\sr{l}{\into}&\law/(l-\psi)\law&\onto&\F_l((\ol T))/\ol\psi(\F_l((T)))\earr}
 with natural vertical maps (which we denote by $\ti{\pht{x}}$). Its
 bottom row is the sequence of cokernels of diagram (1) and thus exact. Its
 right vertical map is an isomorphism, by the definition of $\Xi$ and by
 $\ol\psi(\ol x)=\ol x^l$\,. Consequently, the other
 vertical map, $\Xi\to\law/(l-\psi)\law$, is injective, by the
 snake lemma and $\bigcap_{n\ge 0}l^n\law=0$\,. To finish the proof of the
 lemma we are left with showing the surjectivity of
 $\Xi\to\law/(l-\psi)\law$\,. Starting, in (2), with $\ti
 x\in\law/(l-\psi)\law$ (the middle term in the bottom row) we
 find elements $y_0\in\Xi$ and $\ti x_1\in\law/(l-\psi)\law$ such that $\ti x-\ti y_0=l\ti
 x_1$\,. Continuing, we get $\ti x=\ti y_0+l\ti y_1+l^2\ti
 y_2+\cdots$\,, with $ y_0+l y_1+l^2
 y_2+\cdots\in\Xi$\,.
 \corollary \ $T^2\La=\Xi_2\op(l-\psi)T\La$\Stop
 Since $\Xi\cap\law^+=\Xi_1$, Lemma 2 gives
 $\law^+=\Xi_1\op(l-\psi)\law^+$, i.e.,
 $T\La=\Xi_1\op(l-\psi)T\La$\,. We intersect with $T^2\La$ and
 obtain the corollary from $(l-\psi)T\La\subset T^2\La$ and $\Xi_1\cap
 T^2\La=\Xi_2$\,.
 \Proposition{A} \
 $\lawm=T^\Z\times\mu_{l-1}\times(1+T)^\zl\times\E(\Xi_2)\times(1+l\law)$
 \Stop Given
 $e=\sum_{k\in\Z}e_kT^k\in\lawm$, we will modify $e$ by factors in $T^\Z,\,\mu_{l-1}\times(1+l\law)$ and $(1+T)^\zl$
 to arrive at a new unit $\E(y)$ for some $y\in\Xi_2$. This confirms the claimed product decomposition
 of $\lawm$ but not yet that it is a direct product.
 \ben\item Going modulo $l$, let $
 \ol e=\sum_{k\ge k_0}\ol e_k\ol T^k\in\F_l((\ol T))$ have
 coefficient $\ol e_{k_0}\neq0$. Multiplying $e$ by $T^{-k_0}\in
 T^\Z$ gives a new unit with zero coefficient not divisible by
 $l$ but all coefficients with negative index divisible by $l$; we
 denote it again by $e$.
 \item Now
 $e_0\in\zl\mal=\mu_{l-1}\times(1+l\zl)\subset\mu_{l-1}\times(1+l\law)$, and multiplying $e$ by
 $e_0\me$ allows us to assume that $e=le^-+1+e^+$\,, where $e^-\in\law^-$ and $e^+\in\law^+$\,, so
 $1+e^+\in\La\mal\le\lawm$ and $e(1+e^+)\me=1+l(e^-(1+e^+)\me)\in1+l\law$\,, i.e., $e\equiv1+e^+\mod 1+l\law$.
 \item If
 $1+e^+=1+e_1T+e_2T^2+\cdots$, then multiplying $1+e^+$ by
 $(1+T)^{-e_1}\in(1+T)^\zl$ produces $1+T^2\ti y$ with
 $\ti y\in\La$ (note $(1+T)^z\equiv 1+zT\mod T^2\La$). Hence, by Lemma 1, modulo
 $T^\Z\cdot\mu_{l-1}\cdot(1+T)^\zl\cdot(1+l\law)$, the original unit $e$
 satisfies $e\equiv \E(y')$ with $y'\in T^2\La$.
 \item As $\E(T^2\La)=\E(\Xi_2)\times \E((l-\psi)T\La)$ by the above
 corollary,
 multiplying $\E(y')$ with $\E(y)$ for a suitable $y\in\Xi_2$
 yields an element $\E((l-\psi)y'')$ with $y''\in y_1''T+T^2\La$.
 It follows from Corollary 1 to Lemma 1 that
 $\E((l-\psi)y'')=(1+T)^{-ly_1''}\exp(ly'')$. The first factor is
 in $(1+T)^\zl$ and the second in $1+l\La$.
 \een
 \sn We now prove that we actually have a direct product.
 \sn We have already used $(1+T)^z\equiv1+zT\mod T^2\La$. Together with
 $\E(\Xi_2)\subset 1+T^2\law$ it implies
 that the product $T^\Z\cdot\mu_{l-1}\cdot(1+T)^\zl\cdot\E(\Xi_2)$ is
 direct. Moreover, an element in it which also lies in $1+l\law$
 must equal $\E(y)$ with $y\in \Xi_2$. Indeed, $(1+T)^z\equiv
 1\mod l$ gives $z\equiv0\mod l$, hence $(1+T)^\fra zl\equiv1\mod
 l$, since modulo $l$ we are in characteristic $l$. Thus $z=0$.
 \sn  So assume $\E(y)=1+lz$. Applying
 $\Ll$ gives $y=\Ll(1+lz)=(l-\psi)z'\in(l-\psi)\law$, by Corollary
 1 to Lemma 1. As $y\in\Xi_2\subset T^2\La$, the zero coefficient
 of $z'$ vanishes and the last corollary implies $y=0$.
 This completes the proof of the proposition.
 \definition $\xi:\law=\Xi\op(l-\psi)\law\to\Xi$ is the identity on $\Xi$ and zero
 on $(l-\psi)\law$\Stopp
 \corollary \ We have an exact sequence $$\mu_{l-1}\times(1+T)^\zl\into\law\mal\,\sr{\Ll}{\to}\,\law\onto
 \Xi\,/\,(\Z\cdot\xi((\Ll(T))\op\Xi_2)\,.$$\Stop
 For the proof note that $\xi(\Ll(T))$ is in $\law^-$ and
 non-zero\,: writing $\fra{T^l}{\psi(T)}=\fra{1}{1-lv}$ with
 $v=-\fra1l\sum_{i=1}^{l-1}{l\choose i}T^{-i}$ we have
 $$\barr{l}\Ll(T)=-\log(1-lv)=\sum_{j\ge1}\fra{l^{j-1}}{j}v^j\in\law^-\quad
 \mr{with}\\
 \xi(\Ll(T))\equiv\xi(v)=v\equiv\sum_{i=1}^{l-1}\fra{(-1)^i}{i}T^{-i}\mod
 l\,.\earr$$
 Recall that $\mu_{l-1}\times(1+T)^\zl\subset\ker(\Ll)$, that
 $\Ll\E$ is the identity on $\Xi_2$, and that
 $1+l\law=\exp(l\law)$.
 \sn Suppose now that $e=T^b\ze(1+T)^z\E(x)\exp(ly)$ is in
 $\ker(\Ll)$ (with $b\in\Z,\,\ze\in\mu_{l-1},\,z\in\zl,\,x\in\Xi_2,\,y\in\law$). Then $-b\Ll(T)=x+(l-\psi)y$ implies
 $-b\xi(\Ll(T))=x$ is in $\law^-\cap\Xi_2=0$, hence $b=0=x$ and
 then $y=0$ by 1.~of Lemma 2, as required.
 \sn Concerning $\coker(\Ll)$, it suffices to show that
 $\im(\Ll)=\Z\cdot\xi(\Ll(T))\op\Xi_2\op(l-\psi)\law$. By
 Proposition A, 1.~of Lemma 2 and $\Ll(T)-\xi(\Ll(T))\in(l-\psi)\law$
 this again follows from $\xi(\Ll(T))\notin\Xi_2$.
 \sn This finishes the proof of the corollary.
 \remark When $l=2$, more effort is needed, since
 $-1\in1+2\law$ and `log\,,\,exp' are no longer
 inverse to each other.
 \Section{3}{Kernel and cokernel of $1-\psi$ on $\law$}
 \vspace*{-.75cm}
 \Lemma{3} There is an exact sequence \quad
 $0\to\zl\to\law\sr{1-\psi}{\lto}\law\to(\Xi/\Xi_1)\op\zl\to0$\,.\Stop
 \sn\bmp{5.5cm} We start its proof from the obvious diagram at
 right and show that $\ker(\up)=\F_l$, the constants in $\F_l((\ol
 T))=\law/l\law$. Indeed, \emp\hsp{10}\bmp{9cm}$\barr{ccccc}\law&\sr{l}{\into}&\law&\onto&\law/l\law\\
 \daz{1-\psi}&&\daz{1-\psi}&&\daz{\up}\\
 \law&\sr{l}{\into}&\law&\onto&\law/l\law\earr$ \emp
 $$(\up)(\sum_{k\ge-n}\ol z_k\ol
 T^k)=0\iff\sum_{k\ge-n}\ol z_k\ol T^k=\sum_{k\ge-n}\ol z_k\ol
 T^{lk}=(\sum_{k\ge-n}\ol z_k\ol T^k)^l\,,$$ and the only
 $l-1^\mr{st}$ roots of unity in the field $\F_l((\ol T))$ are the
 constants $\neq0$. The above implies $\ker(1-\psi)=\zl+l\ker(1-\psi)$. By
 successive approximation this gives $\ker(1-\psi)=\zl$\,.
 \sn Turning back to the diagram, we obtain from the snake
 lemma the short exact sequence
 $$\coker(1-\psi)\sr{l}{\into}\coker(1-\psi)\onto\coker(\up)\,.$$ We
 compute its right end. Because $\F_l((\ol T))$ is complete in the
 $v_{\ol T}$-topology, $\sum_{n\ge0}\ol z^{l^n}$ converges for
 every element $\ol z=\sum_{k\ge1}\ol z_k\ol T^k$, hence
 $(\up)(\sum_{n\ge0}\ol z^{l^n})=\ol z$ implies that these $\ol
 z$ all belong to $\im(\up)$. Also, $\ol T^i-\ol T^{li}=(\up)(\ol
 T^i)\in\im(\up)$. Thus, $\coker(\up)$ is spanned by the images of
 $\ol T^j$ with $j=0$ or $j<0\ \&\ l\nmid j$. These elements are
 actually linearly independent over $\F_l$. To see this, read an
 equation $$\sum_{{-n\le k<0}\atop{l\nmid k}}\ol z_k\ol
 T^k+z_0=(\up)(\ol x)=\sum_{-n\le k<0}\ol x_k(\ol T^k-\ol T^{lk})$$
 coefficientwise from $k=-n$ to $k=0$.
 \sn Going back to the short exact sequence displayed above, we
 now
 realize that $\Xi/\Xi_1\,\op\zl$ maps onto $\coker(1-\psi)$, since $\law$
 is $l$-complete. And by the last paragraph, this surjection
 is, in fact, an isomorphism.
 \Section{4}{Kernel and cokernel of $\Ll$ on $\lawh$}
 As in the introduction, $H$ is a finite abelian $l$-group and
 $\lawh$ is its group ring over $\law$. Perhaps the description \ $\lawh=\zl[[\Ga\times H]]_\wedge$\,, with $\Ga$
 denoting the cyclic pro-$l$ group generated by $1+T$, gives a better
 understanding of the ring homomorphism $\psi$ on $\lawh$\,: $\psi$ is
 induced by \ $\psi(g)=g^l$ for $g\in\Ga\times H$\,. And the
 integral logarithm \ $\Ll:\lawh\mal\to\lawh$\,, as before, takes a unit
 $e\in\lawh\mal$ to \
 $\Ll(e)=\fra1l\log\fra{e^l}{\psi(e)}$\,.
 \sn For the discussion of its kernel and cokernel we first invoke
 the augmentation map \ $\lawh\to\law\,,\,h\mapsto1$ for $h\in
 H$\,, so that we can employ our earlier results. Let $\fg$ denote
 its kernel and note that $1+\fg\subset\lawh\mal$, as
 $\fg\subset\fr\df\mr{rad}(\lawh)=\fg+l\lawh$\,; moreover, for the same reason, $\lawh\mal\to\lawm$ is
 surjective.
 \Proposition{B} \ $\Ll:\lawh\mal\to\lawh$ has
 \ben\item[]$\ker(\Ll)=\mu_{l-1}\times(1+T)^\zl\times H
 \quad(\,=\,\mu_{l-1}\times(\Ga\times H)\,)$\,,
 \item[] and $\coker(\Ll)$ is described by the split exact
 sequence $$(\Xi/\Xi_1\op\zl)\ot_\zl
 H\into\coker(\Ll)\onto\Xi\,/\,(\Z\xi(\Ll(T))\op\Xi_2)\ .$$\een\Stop
 The proof begins with the commutative diagram
 $$\barr{ccccc}1+\fg&\into&\lawh\mal&\onto&\lawm\\
 \daz{\Ll}&&\daz{\Ll}&&\daz{\Ll}\\
 \fg&\into&\lawh&\onto&\law\earr$$ with exact rows which are
 split by the {\em same} inclusion $\law\into\lawh$ of rings.
 Here the right square commutes because $\psi$ and `log' both
 commute with augmentation, and thus induces the left square since
 the sequences are exact.
 \enlargethispage{5mm}
 \sn The right vertical $\Ll$ fits into the exact sequence of
 the corollary to Proposition A. Similarly we will need
 \Lemma{4} There is an exact sequence \quad
 $H\into1+\fg\sr{\Ll}{\lto}\fg\onto(\Xi/\Xi_1\op\zl)\ot_\zl H\,.$
 \Stop
 Proposition B follows from Lemma 4 and the snake lemma\,: for
 $\mu_{l-1}\times(1+T)^\zl\times H\subset\ker(\Ll)$ maps onto the
 kernel of the right vertical $\Ll$\,, and the cokernel sequence
 splits because the natural splittings in the commutative diagram
 are compatible. So it remains to prove Lemma 4, which we do next.
 \ben\item[a)] \  {\em $\fg/\fg^2\simeq\law\ot_\zl H$
 by \ $h-1\mod\fg^2 \ \mapsto\, h$}
 \nz This is a consequence of
 $\lawh=\law\ot_\zl\zl[H]$ and the natural isomorphism $\De
 H/\De^2H\simeq H\,,\,h-1\mapsto h$, where
 $\De H=\group{h-1:h\in H}_\zl$ is the augmentation ideal of the
 group ring $\zl[H]$, so $\fg=\law\ot_\zl\De H$.
 \item[b)] {\em If $e=1+\sum_{1\neq h\in H}e_h(h-1)\in1+\fg$ (with $e_h\in\law$),
 then
 \vspace*{-3mm}
 \bct$\Ll(e)\equiv \sum_h(e_h-\psi(e_h))(h-1)\mod\fg^2\,.$\ect}
 \vspace*{-3mm}
 Indeed, modulo $\l\fg^2$ we have
 {\small $$\barr{l}e^l\equiv
 1+l\sum_he_h(h-1)+\sum_he_h^l(h-1)^l\equiv1+l\sum_he_h(h-1)+\sum_h\psi(e_h)(h-1)^l\equiv\\
 \,\,1+l\sum_he_h(h-1)\!+\!\sum_h\psi(e_h)(h^l-1)\!-\!l\sum_h\psi(e_h)(h-1)\equiv\psi(e)\!+\!l\sum_h(e_h-\psi(e_h))(h-1)\,,\earr$$}
 \hsp{-1.2}so
 $\fra{e^l}{\psi(e)}\equiv1+\psi(e)\me l\sum_h(e_h-\psi(e_h))(h-1)\equiv
 1+l\sum_h(e_h-\psi(e_h))(h-1)\mod l\fg^2$ as
 $\psi(e)\me\in1+\fg$. Now apply\, `$\!\fra1l\log$\,'.\een
 From a),\,b) we get the right square of the commutative diagram
 $$\barr{ccccc}1+\fg^2&\into&1+\fg&\onto&\law\ot_\zl H\\
 \daz{\Ll}&&\daz{\Ll}&&\daz{(1-\psi)\ot1}\\
 \fg^2&\into&\fg&\onto&\law\ot_\zl H\earr$$
 with left square induced by the exactness of the rows. The map
 $(1-\psi)\ot1$ has kernel and cokernel given by tensoring the sequence in Lemma 3
 with $H$\,: it remains exact since it is composed of two short
 exact sequences of torsionfree $\zl$-modules. So the snake lemma
 reduces Lemma 4 to proving that $\Ll:1+\fg^2\to\fg^2$ is an
 isomorphism.
 \sn We do this by induction on
 $|H|$ and, to that
 end, choose an element $h_0\in H$ of order $l$ and let $H\,\ti\to\,\ti
 H=H/\group{h_0}$ be the natural map.
 \sn Recalling that $\fr=\mr{rad}(\lawh)=\fg+l\lawh$\,,
 we start with the right square of the diagram
 $$\barr{cccccc}1+(h_0-1)\fr&\into&1+\fg^2&\onto&1+\ti\fg^2&\\
 \daz{\Ll}&&\daz{\Ll}&&\daz{\ti\Ll}&\\
 (h_0-1)\fr&\into&\fg^2&\onto&\ti\fg^2&,\earr$$
 which commutes since $\psi$\,,\,`log' commute with\, \~\ . Since
 $\ti\Ll$ is an isomorphism by the induction hypothesis, it
 suffices to show that the kernels in the rows are as shown and
 that the left $\Ll$ is an isomorphism\,:
 \ben
 \item[i.] {\em $\fg^2\to\ti\fg^2$ has kernel $(h_0-1)\fr$}\,. Since
 $(h_0-1)\fr$ is in the kernel of\, \~\  and $l(h_0-1)$ is in
 $\fg^2$, by $l(h_0-1)\equiv h_0^l-1\mod\fg^2$, it remains to
 check
 $$(h_0-1)\lawh\cap\fg^2\subset(h_0-1)\fr\,.$$
 If $(h_0-1)b=(h_0-1)\sum_{h\in H}b_hh\in\fg^2$ (with
 $b_h\in\law$), then the isomorphism $\fg/\fg^2\simeq\law\ot_\zl H$ takes
 $(h_0-1)b$ to $0=\sum_{h\in H}b_h\ot h_0=(\sum_{h\in H}b_h)\ot
 h_0$\,, whence $\sum_{h\in H}b_h\in l\law$, since $h_0$ has order
 $l$. Thus, $(h_0-1)b\in(h_0-1)(\sum_{h\in
 H}b_h(h-1)+l\law)\subset(h_0-1)\fr$\,.
 \een The same argument applies to the kernel in the top row. It
 follows that $\Ll(1+(h_0-1)\fr)\subset(h_0-1)\fr$\,.
 \ben
 \item[ii.] {\em $\Ll:1+(h_0-1)\fr\to(h_0-1)\fr$ is an isomorphism}\,. If $x\in\fr$, then
 $\psi(h_0-1)=0$ implies that $\Ll(\exp((h_0-1)x))=(h_0-1)x$\,, hence $\Ll$ is
 onto, and \
 $\Ll(1-(h_0-1)x)=$ $$\log(1-(h_0-1)x)=-(h_0-1)(x-x^l)+(h_0-1)^2x^2\la_x=-(h_0-1)x+(h_0-1)x^2\la_x'$$
 with some $\la_x,\la_x'\in\lawh$ by ($\dag$) in [RW3, p.40] (with $z$
 replaced by $h_0$). If this is zero, then
 $(h_0-1)x(1-x\la_x')=0$ with $1-x\la_x'\in\lawh\mal$.
 So $\Ll$ is injective.\een
 \remark Admittedly, Proposition B is closer to the corollary to
 Proposition A than to Proposition A itself, as $\lawh\mal$ has not
 been determined.

 \enlargethispage{.75cm}
 \bbn{\large {\sc References}}
 \footnotesize
 \bn
 \btb{rp{13cm}}
 \,[CR]   & Curtis, C.W.~and Reiner, I., {\em Methods of
            Representation Theory,
            I,II.} John Wiley \& Sons (1981,1987) \\
 \,[C]    & Coleman, R., {\em Local units modulo circular units.} Proc.\,AMS {\bf 89} (1983), 1-7\\
 \,[FeKu] & Fesenko, I., Kurihara, M., {\em Invitation to higher local fields.} Geometry \&\ Topology Monographs
            {\bf 3} (2000), ISSN 1464-8997 (on line)\\
 \,[F]    & Fröhlich, A., {\em Galois Module Structure of Algebraic
            Integers.} Springer-Verlag (1983)\\
 \,[FK]   & Fukaya, T., Kato, K., {\em A formulation of conjectures on $p$-adic zeta
            functions in non-commutative Iwasawa theory.}
            Proc.\,St.\,Petersburg Math.\,Soc.\,{\bf 11} (2005)\\
 \,[K]    & Kato, K., {\em Iwasawa theory of totally real
            fields for Galois extensions of Heisenberg type.} Preprint (`Very preliminary version'\,, 2006)\\
 \,[L]    & Lang, S., {\em Cylotomic Fields I-II.} Springer GTM
            {\bf 121} (1990)\\
 \,[O]    & Oliver, R., {\em Whitehead Groups of Finite Groups.}
            LMS Lecture Notes Series {\bf 132}, Cambridge (1988)\\
 \,[RW2,3]  & Ritter, J.~and Weiss, A., {\em Toward equivariant Iwasawa theory, II; III.}
           Indagationes Mathematicae {\bf 15} (2004), 549-572; Mathematische Annalen {\bf 336}
           (2006), 27-49  \\
\,[RW4]    & ------------------------\,, {\em Non-abelian
           pseudomeasures
           and congruences between abelian Iwasawa $L$-functions.}
           To appear in Pure and Applied Mathematics Quarterly (2007)\\
\,[RW5]   & ------------------------\,, {\em Congruences between
           abelian pseudomeasures.}
           Preprint (2007)\\
\,[RW6]   & ------------------------\,, {\em Equivariant Iwasawa
           theory\,: an example.}
           Preprint (2007)\\
 \etb

 \bn {\footnotesize \bct Institut für Mathematik $\cdot$
 Universität Augsburg $\cdot$ 86135 Augsburg $\cdot$ Germany \\
 Department of Mathematics $\cdot$ University of Alberta $\cdot$
 Edmonton, AB $\cdot$ Canada T6G 2G1   \ect

\end{document}